\documentclass[a4paper]{article}
\usepackage{latexsym,amsfonts,amscd,amssymb,enumerate,amstext,amsmath,dsfont}
\usepackage[sc,tiny,center,compact]{titlesec}
\usepackage{amsthm}
\usepackage[all]{xy}
\usepackage{titletoc}
\titlecontents{section}[0pt]{\small\scshape}{\contentslabel{1.5em}}{}{\titlerule*[0.75pc]{.}\contentspage}[\textbullet\ ]

\titlecontents{subsection}[1.5em]{\small\scshape}{\contentslabel{2.3em}}{}{\titlerule*[0.75pc]{.}\contentspage}[\textbullet\ ]
\titlecontents{subsubsection}[3.8em]{\small\scshape}{\contentslabel{3.0em}}{}{\titlerule*[0.75pc]{.}\contentspage}[\textbullet\ ]

\author{M. Grime \\ University of Bristol\\ email: Matt.Grime@bris.ac.uk}
\title{The correct relatively stable category for idempotent modules\footnote{2000MSC: 18E05,18E20,18E30,18E35,18G25,20C20,20J05. Last updated 13/10/2007}}
\date{}
\newtheoremstyle{ordinary}{1ex}{0pt}{}{}{\scshape}{.}{\newline}{}
\theoremstyle{ordinary}
\newtheorem{thm}{Theorem}[section]
\newtheorem{defn}[thm]{Definition}
\newtheorem{lem}[thm]{Lemma}
\newtheorem{cor}[thm]{Corollary}

\newtheorem{prop}[thm]{Proposition}

\newcommand{\catt}{\mathcal{T}}
\newcommand{\catl}{\mathcal{L}}
\newcommand{\cats}{\mathcal{S}}
\newcommand{\catc}{\mathcal{L}'}
\renewcommand{\mod}{\mathrm{mod}(kG)}
\newcommand{\Mod}{\mathrm{Mod}(kG)}
\newcommand{\stmod}{\mathrm{stmod}_w(kG)}
\newcommand{\Stmod}{\mathrm{StMod}_w(kG)}

\begin{document}
\maketitle

\begin{abstract}
We answer a question posed in \cite{virtual}, and demonstrate that in general Rickard modules in relatively stable categories are not idempotent modules even if one localizes with respect to a tensor ideal subcategory. We also show that there is a modification one can make so as to recover the idempotent behaviour.

\end{abstract}

\section{Introduction}

It is now more than a decade since Rickard introduced his idempotent modules in \cite{idempotent}. They arise from performing a Bousfield localization in $\mathrm{StMod}(kG)$.

\begin{prop}\label{idem}Suppose that $\catc$ is a thick subcategory of $\mathrm{stmod}(kG)$. Let $\catl$ be the smallest localizing subcategory of $\mathrm{StMod}(kG)$ containing $\catc$, then for each $M$ in $\mathrm{StMod}(kG)$ there is a so-called Bousfield triangle
\[ T_\catl(M):=  M_\catl \to M \to M_{\catl^\perp} \leadsto \]
in $\mathrm{StMod}(kG)$ with $M_\catl$ in $\catl$ and  $(\catl,M_{\catl^\perp})=0$. Moreover, if $\catc$ is a tensor closed subcategory (or smashing in the language of topologists) of $\mathrm{stmod}(kG)$, then three further conditions are satisfied:

\begin{enumerate}[(i)] 
\item for any object $M$ the  triangle $T_\catl(M)$  is isomorphic to $M\otimes T_\catl(k)$;
\item the object $M_\catl \otimes_k M_{\catl^\perp}$ is zero in $\mathrm{StMod}(kG)$;
\item in $\mathrm{StMod}(kG)$, $M_\catl\otimes_k M_{\catl}\cong M_\catl$.
\end{enumerate}
\end{prop}
We remind the reader that  $\catl'$ is tensor closed if $x\in \catl'$ and $y \in \mathrm{stmod}(kG)$ implies that $x\otimes y$ is in $\catl'$. It is the last of these properties that justifies the name idempotent module. We will refer to these three extra conditions as \emph{tensor ideal property} (or TIP for brevity) (i), (ii), and (iii) of the localization. 

It is natural to ask if Bousfield localizations can be performed in any of the relatively stable triangulated categories one may associate to $kG$. Indeed, the fact that one can construct Bousfield triangles in relatively stable categories appears in \cite{virtual}[Prop.~6.5] without proof. In that paper they comment that it is not clear whether the triangles so constructed satisfy TIP~(i) and suggest that this should be the object of further study.

We demonstrate two things. First, even if $\catl'\subseteq \stmod$ is tensor closed, then TIP~(i) may fail if we localize in $\Stmod$. Second, if one chooses a different category in which to localize, then one can create a theory of idempotent modules for tensor closed subcategories for which we do recover TIP~(i)--(iii).

This relativized statement of proposition \ref{idem} needs careful formulation, more so than the ordinary stable category. In some sense this because $\Stmod$ is the wrong object to study. Given a compactly generated localizing subcategory, $\catl\subseteq\Stmod$, one can invoke Bousfield's theorem. This case is frequently called \emph{finite} Bousfield localization. If $\catl$ is the smallest localizing category containing $\catl' \subseteq \stmod$, then we are certainly in such a situation, and the Bousfield triangles $T_\catl(M)$ will exist. In the ordinary stable case one can be laxer in the statement of the existence of Bousfield triangles: since the category $\mathrm{StMod}(kG)$ is itself compactly generated, one can Bousfield localize with respect to any localizing subcategory.  This result does not generalize: we are restricted to considering only finite Bousfield localizations where $\catl$ is compactly generated.

 We will demonstrate this failure of  TIP~(i)  by example. The example indicates a necessary condition for when one has a relativized statement of proposition \ref{idem}. We will also prove that this is sufficient.

We start with a reminder to the reader of the technical conditions needed so that one may invoke Bousfield's localization argument. We discuss in section \ref{examples} which aspects of these technical conditions are known to be true in the relatively stable triangulated categories one may associate to $kG$ . We offer two examples which make it clear that we need to be careful about the kinds of localizing subcategories we can use if we wish to work in $\Stmod$.

In section \ref{locbig} we use the results of \ref{examples} to show that TIP~(i) can and does fail.  In section \ref{locsmall} we will argue that one should think of $\stmod^\oplus$ as the correct triangulated category to work with if one wishes to correctly generalize Rickard's construction. Indeed we show that one can replace the ordinary stable category with the relatively stable category in proposition \ref{idem} if and only if 

\[ \Stmod \cong \stmod^\oplus,\]
and in the process we  will demonstrate that this is a non-trivial fact to verify.

\section{Bousfield localization}

We assume that $\catt$ is a triangulated category with arbitrary direct sums.  Throughout, we will assume that all subcategories are closed under isomorphism. We will use the notation $(X,Y)_\catt$ for the Hom sets in $\catt$, omitting the subscript whenever we may do so without fear of confusion. A full subcategory $\catl$ is defined to be \emph{localizing} if it is closed under taking $\catt$-direct sums, and \emph{thick} if it is closed under direct summands.  If $\cats$ is any subcategory we can define two localizing subcategories, $\cats^\oplus$ and $\cats^\perp$. $\cats^\oplus$ is the smallest localizing subcategory of $\catt$ containing $\cats$, and $\cats^\perp$ is the class of objects that satisfy
\[\cats^\perp:=\{ X \in \catt : (\cats,X)_\catt =0\}.\]

Recall that an object $c$ in $\catt$ is compact  or small (thus we use lower case letters to distinguish compact objects) if the canonical map
\[  \mathop{\oplus}_{\lambda \in \Lambda} (c,X_\lambda)_\catt  \to (c,\mathop{\oplus}_{\lambda \in \Lambda} X_\lambda)_\catt\]
is an isomorphism. Or, more descriptively, $c$ is compact if any map to a direct sum factors through a finite subsum.
If $\cats$ is a class of objects in $\catt$, then we say it generates $\catt$ if 
\[(S[n],X)_\catt = 0\  \forall \ S \in \cats \ \mathrm{and}\ n \in \mathbb{Z} \iff X=0.\]
Piecing these two definitions together we have the definition of compact generation: a triangulated category $\catt$ is compactly generated if there is a \emph{set} of compact objects that generate. We can now state Bousfield's theorem.

\begin{thm}Suppose that $\catt$ is a triangulated category with direct sums and that  $\catl$ is a localizing subcategory of $\catt$. Further suppose  that one of $\catl$ or $\catt$ is compactly generated, then, \emph{when it exists}\footnote{There are some issues in the case that $\catt$ is a compactly generated over whether the Hom sets of the quotient exist in this universe. We do not address these issues, since we will only care about the case when $\catl$ is compactly generated when there are no issues}, the Verdier quotient functor $Q:\catt \to \catt/\catl$ has a left adjoint. Moreover  for each object $X$ in $\catt$ there is a distinguished triangle

\[ X_\catl \to X \to X_{\catl^\perp} \leadsto \]
and $X_\catl$ is in $\catl$, and $X_{\catl^\perp}$ is in $\catl^\perp$. The map $X_\catl \to X$ is given by the counit of the adjunction.
\end{thm}

An instant corollary of this is that every map from an object in $\catl$ to $X$ factors uniquely through $X_\catl$, and that this characterizes the triangle up to isomorphism.

\section{Relatively stable module categories}\label{examples}

Let $k$ be an algebraically closed field of charactersitic $p$ which divides the order of the (finite) group $G$. We shall consider triangulated quotients of the categories $\mod$ and $\Mod$ of finite dimensional and arbitrary dimensional $kG$ modules. We will use lower case letters for finite dimensional modules. It will be  clear that  they are  compact  in the triangulated quotients. We will define our relatively stable categories in a slightly novel way, though these relatively stable categories should now be considered well known in the representation theoretic world.

\begin{defn} Let $w$ be an element of $\mod$. Define a class of objects in $\Mod$ as follows

\[ X \in \mathcal{P}(w)  \iff  \exists Y \in \Mod \ \mathrm{s.t.}\ X | w\otimes Y.\]
\end{defn}
The class $\mathcal{P}(w)$ is pre-enveloping, pre-covering, contains the ordinary projective modules and passing to the stable category there is a quotient of \emph{additive} categories

\[  \Stmod:=\mathrm{StMod}(kG)/\mathcal{P}(w)\]
where the objects are the same and  the Hom sets are the quotients of those in $\mathrm{StMod}(kG)$ by those that factor through an object in $\mathcal{P}(w)$. We emphasise that this na\"{\i}ve construction is not the Verdier quotient. However, the quotient is triangulated. The shift functor is given by completing the morphism $P(X)\to X$ to a triangle
\[ \Omega_w(X)\to P(X) \to X \leadsto \]
where $P(X)$ is a right $\mathcal{P}(w)$-approximation to $X$. The distinguished triangles in $\Stmod$ are the triangles in $\mathrm{StMod}(kG)$ that split on tensoring with $w$.

 Note that we can replace $\mathrm{Stmod}(kG)$ with $\mathrm{stmod}(kG)$.
We will use $\stmod$ and $\Stmod$ for the triangulated quotients. Notice that $\stmod$ is a full triangulated subcategory of $\Stmod$. 

\medskip

\noindent
\textsc{Remarks}
\begin{enumerate}
\item  We can recover the ordinary stable categories by setting $w=kG$, and the normal notion of relatively projective with respect to a subgroup by choosing $w= \mathrm{Ind}_H^G(k)$.

\item These categories are trivial, in the sense that every object is $w$-projective, if $w$ has any summand of dimension prime to $p$.

\item For consistency with the work of others, and clarity, we will make explicit the $w$ dependence in our discussions. Thus we will talk of the relatively stable category (for an arbitrary $w$) and reserve the term ordinary stable category for the case $w=kG$. The reader should never have a situation where she finds herself wondering what we mean by `stable category.'

\item We choose to give this definition (which is easily seen to be equivalent to the normal definition as a quotient of the module categories) since the author has been asked  many times now what the relation between $\stmod$ and $\mathrm{stmod}(kG)$ is: clearly the quotient from $\mod$ to $\stmod$ must factor through the ordinary stable category, but this  factorization is only at the level of additive functors. 

\end{enumerate}
Bousfield localization for ordinary stable module categories first appears in \cite{idempotent}, and has since proven to be a very powerful tool in the representation theorist's armoury: \cite{classify}, \cite{direct} are but two examples. The latter includes more references. There are several subtle, and not so subtle, differences between the ordinary and the relative case.  For instance, the direct limit in $\Mod$ of flat modules is a flat module. Thus the direct limit of ordinary projective modules is an ordinary projective module since flat and projective are synonymous for group algebras (of finite groups). This fails in the relative case, and such differences lead to interesting new phenomena to investigate.

In order to perform a finite Bousfield localization we need to know that there are compact objects to generate localizing subcategories. It is trivial to see that there are many compact objects for us to use:
\begin{prop}If $c$ is isomorphic in $\Stmod$ to a finite dimensional object then it is compact.
\end{prop}

\medskip
If we wished to localize with respect to arbitrary localizing subcategories, then we would need to know that the relatively stable module category itself is compactly generated. There are several  issues to consider.

\begin{enumerate}[(i)]
\item Do  the simple objects form a compact generating set?
\item If not is there a set of finite dimensional objects that are generators?
 \item Are there more compact objects which we can consider, if necessary?
\end{enumerate}
\noindent
Only partial answers to these questions are known. In order they are

\begin{enumerate}[(i)]
\item In general, no, they are not sufficient:  suppose that $s$ is a simple module with relatively injective hull $I(s)$, and that the socle of $I$ is $s$, let $s'$ be some simple module such that there is  a non-split extension 
\[0\to I(s) \to X \to s'\to 0\] 
so that $s$ is the socle of $X$. There are no non-zero morphisms in the relatively stable category from any simple object into $X$: the only simple with a non-zero map to $X$ in the module category is $s$; the inclusion factors through $I(s)$. In the example of theorem \ref{compact} below, $\Omega(k)$ is  such an $X$.
\item Such a generating set is known to exist in many cases (recent work of the author with Peter Jorgensen) when $\mathcal{P}(w)$ satisfies some appropriate finiteness condition. We  offer an example below. We also offer a counter example.
\item This is completely open, as far as the author is aware, except for those cases where one can apply Thomason localization in the previous answer. A constructive proof that the finite dimensional objects are precisely, up to isomorphism, the compacts does not exist. The author knows of no non-vacuous cases of infinite dimensional module that is compact in $\Stmod$.
\end{enumerate}

\noindent
It is the first two of these topics that we are interested in in this note. The next theorem illustrates that we can have compact generation.

\begin{thm}\label{compact}
Let $G=C_2\times C_2$, $H=C_2\times e$, and $w=\mathrm{Ind}_H^G(k)$, then $\Stmod$ is compactly generated. In fact $\Omega(k)$ (the ordinary Heller translate) is a compact generator.
\end{thm}
\begin{proof}
If the reader is unfamiliar with the module theory of $G$ we direct her attention to \cite{direct}. There are only two relatively projective indecomposables. The $4$ dimensional free module and a uniserial $2$ dimensional module. One can verify that any indecomposable $kG$-module with at least two copies of $k$ in its socle is detected\footnote{We say $(X,?)_\catt$ detects $Y$ precisely when $(X,Y)_\catt$ does not vanish} by $(k,?)$. There are only five $kG$-modules with one copy of $k$ in the socle, so clearly there is a set of generators. We can do even better though: of these, two are relatively projective, two are detected by $(k,?)$, and one, $\Omega(k)$ is not. Thus $\{k,\Omega(k)\}$ is a set of compact generators, and it is easy to see that one does not need $k$ as there is a morphism $\Omega(k)\to k$ that is not zero in $\stmod$.
\end{proof} 

Recent work with Peter Jorgensen has extended this result to a much larger class of examples, and will appear at some future date.

However, compact generation by the finite dimensional objects should be seen as the exception rather than the norm. It would imply that the class $\mathcal{P}(w)$, considered in the module category, is closed under arbitrary direct limits, as we will see in the next theorem.  Such a class is called \emph{definable}. The reader is referred to \cite{krause}[Ch.2] for a discussion on  definable classes.

The next example is more indicative of the kinds of behaviour one expects. It shows that there are very simple examples where the finite dimensional objects do not form a generating set.

\begin{thm}\label{noncomp} Let $H\cong K \cong L \cong C_p$, and set $G=H\times K \times L$. Further, let $w=\mathrm{Ind}_{H\times K}^G(k)$. Then $\Stmod$ is not generated by the finite dimensional objects.\end{thm}
\begin{proof}We only sketch an idea of the proof. One may find a sequence of relatively projective modules whose direct limit  $M$ is not relatively projective. Any map from a finite dimensional object to $M$ must factor through one of the finite dimensional relatively projective submodules and is thus zero in $\stmod$.
\end{proof}
The proof is not hard, but it is lengthy. For more details, the reader is referred to \cite{noncompact},  which originally appears as an appendix in \cite{thesis}.

\section{Bousfield localization in $\Stmod$}\label{locbig}

From the discussion of the preceding section it follows that we have enough compact objects to apply Bousfield localization in $\Stmod$, but we are generally only able to perform finite Bousfield localizations in $\Stmod$.  The localization will give idempotent \emph{functors} (in the sense that $(M_\catl)_\catl\cong M_\catl$). 

If we start with $\catc \subseteq \stmod$, then $\catl$, the smallest localizing subcategory containing $\catc$, is compactly generated and one obtains Proposition 6.5 of \cite{virtual}. There is no reason to suppose at this stage  that we recover the TIPs when $\catl'$ is tensor closed. In fact, we can go further than this and show that there are cases where the idempotent module behaviour must fail. We apply the following easy proposition to our example \ref{noncomp}.
\begin{prop}
If $\Stmod$ is such that for all (tensor ideal) subcategories $\catl'\subseteq\stmod$ every Bousfield triangle $T_\catl(M)$ is isomorphic to $M\otimes T_\catl(k)$, then $\Stmod\cong \stmod^\oplus$.

\end{prop}
\begin{proof}

Suppose we choose $\catl'$ to be $\stmod$, then $\catl'$ is trivially a tensor closed subcategory of $\stmod$. Since $k$ is in $\catl$, we  must have 
\[T_\catl(k) = k\to k \to 0 \leadsto \]
as the Bousfield triangle associated to the trivial module. If we had the tensor property of triangles, then it follows that all modules $M$ lie in $\stmod^\oplus$. 
\end{proof}

By the example given in theorem \ref{noncomp} we know there exist choices of $G$, $w$, and $\catl'$ such that up to isomorphism there is a proper containment
\[  \stmod^\oplus \subset \Stmod. \]
Thus we have shown a non-vacuous necessary condition that $\Stmod$ must satisfy in order to obtain idempotent modules for all localizations with respect to tensor ideal thick subcategories.  In section \ref{locsmall} we will show that this is also a sufficient condition, thus answering the question posed after \cite{virtual}[Prop~6.5].



\section{Bousfield localizations in $\stmod^\oplus$}\label{locsmall}

Suppose that instead of looking at finite Bousfield localizations in $\Stmod$, we choose to fix $\catt:=\stmod^\oplus$ as the ambient category for the rest of this note. $\catt$ is now compactly generated, and by Thomason localization the compact objects are precisely those isomorphic to direct summands of finite direct sums of finite dimensional modules. With this restriction one obtains a theory much more in the spirit of Rickard's construction.

\begin{prop}
Suppose that $\catc$ is a thick subcategory of $\stmod$. Define $\catl$ to be the smallest localizing subcategory of $\catt$ containing $\catc$. Then there  are Bousfield triangles
\[ T_\catl(M):= M_{\catl}\to M \to M_{\catl^\perp} \leadsto \]
in $\catt$ given by localizing with respect to $\catl$. If $\catc$ satisfies the further condition that it is a tensor closed ($x$ in $\catc$ and $y$ in $\stmod$ implies $x\otimes_k y$ is in $\catc$) then three extra conditions hold:
\begin{enumerate}[(i)] 
\item the Bousfield triangle $T_\catl(M)$ for any module $M$ is isomorphic to $M\otimes T_\catl(k)$;
\item the object $M_\catl\otimes_k M_{\catl^\perp}$ is zero in $\cats$;
\item in $\catt$, $M_\catl\otimes_k M_{\catl}\cong M_\catl$.
\end{enumerate}
\end{prop}
Let us state the important corollary of this proposition before we conclude this note with its proof.
\begin{cor}
The Bousfield triangle $T_\catl(M)$ is isomorphic to $M\otimes T_\catl(k)$ for all $M$ and tensor ideal thick subcategories $\catl'\subseteq\stmod$ if and only if 
\[ \Stmod \cong \stmod^\oplus.\]
\end{cor}
We now prove the proposition.
\begin{proof}Of course, we need not show the existence of $T_\catl(M)$. The proof of TIP~(i), (ii), and (iii)  follows from a series of short lemmas that are essentially formal consequences of the definition of $\catt$, and appear in \cite{idempotent}[Sec.~5]. We are merely changing that which needs to be changed.
\end{proof}

\begin{lem}
Suppose that $m$ is in $\stmod$. If $\catc$ is tensor closed, and $X$ is in $\catl$, then $m\otimes X$ is in $\catl$. 
\end{lem}
\begin{proof} Define a class of objects $\{ Y \!: \! Y\otimes m \in \catl\}$. This contains $\catc$, and is closed under under sums and triangles, hence it contains $\catl$.
\end{proof}
\begin{lem}Suppose that $M$ is in $\catt$. If $\catc$ is tensor closed, and if $X$ is in $\catl$, then $M\otimes X$ is in $\catl$.
\end{lem}
\begin{proof}
The class of objects for which this is true contains $\stmod$, by the previous lemma, and is closed under direct sums and triangles, hence is all of $\catt$.
\end{proof}
Given these arguments, the reader can safely be left to complete the proof by following \cite{idempotent}[Lem.~5.10, 5.12, Prop.~5.11, 5.13].

Note that our final lemma is the key point where any hope of getting idempotent modules in the general case disappears: the triangulated closure of $\stmod$ in $\Stmod$ may well be a proper (inequivalent) subcategory. Of course, this naturally raises more questions about when the inclusion is an equivalence of categories. At this stage, our lack of understanding of when the finite dimensional objects generate $\Stmod$ is the major obstruction. At the moment we have sufficient conditions, but by no means a thorough understanding of necessary ones. An alternate line of attack would be to classify compact objects in $\Stmod$.  Both of these should be the focus of further research.

\section{Acknowledgements}
The author would like to thank Jeremy Rickard, Joseph Chuang, Thorsten Holm and Peter Jorgensen for their patience in listening to his explanations of the ideas herein as they mutated into their current form. 

\begin{small}
\nocite{neeman}
\nocite{virtual}
\nocite{ideals}
\nocite{idempotent}
\nocite{okuyama}
\nocite{bousfield}
\bibliographystyle{plain}
\bibliography{../masterbibliography}

\begin{thebibliography}{10}

\bibitem{classify}
D.~J. Benson, J.~F. Carlson, and J.~Rickard.
\newblock Thick subcategories of the stable module category.
\newblock {\em Fundamentae Mathematica}, (153):pp 59--80, 1997.

\bibitem{direct}
D.~J. Benson and Wayne Wheeler.
\newblock Direct sum decompositions of infinitely generated modules.
\newblock {\em Trans. Amer. Math. Soc.}, 351:pp 3843--3855, 1999.

\bibitem{ideals}
Jon~F. Carlson and Chuang Peng.
\newblock Relative projectivity and ideals in cohomology rings.
\newblock {\em J. of Algebra}, 183:pp 929--948, 1996.

\bibitem{virtual}
Jon~F. Carlson, Chuang Peng, and Wayne Wheeler.
\newblock Transfer maps and virtual projectivity.
\newblock {\em J. Algebra}, 204(1):pp 286--311, 1998.

\bibitem{precover}
Matthew Grime.
\newblock Precovers, localizations and stable homotopy theory.
\newblock arxiv:0708.2866.

\bibitem{thesis}
Matthew Grime.
\newblock {\em Triangulated Categories, Adjoint Functors and Bousfield
  Localization}.
\newblock PhD thesis, University of Bristol, 2005.

\bibitem{noncompact}
Matthew Grime.
\newblock Finite dimensional modules and perpendicular subcategories.
\newblock arxiv:0708.3329, June 2007.

\bibitem{krause}
Henning Krause.
\newblock The spectrum of a module category.
\newblock {\em Memoirs of the American Mathematical Society}, 149:125 pages,
  2001.
\newblock Habilitationsschrift, Bielefeld 1998.

\bibitem{neeman}
Amnon Neeman.
\newblock {\em Triangulated Categories}.
\newblock Princeton University Press, 2001.

\bibitem{okuyama}
Tetsuro Okuyama.
\newblock A generalization of projective covers of modules over finite group
  algebras.
\newblock unpublished manuscript.

\bibitem{idempotent}
J.~Rickard.
\newblock Idempotent modules in the stable category.
\newblock {\em J. London Math. Soc.}, 56(1):pp 149--170, 1997.

\bibitem{bousfield}
J.~Rickard.
\newblock Bousfield localization for representation theorists.
\newblock {\em Trends in Math., Birkh\"{a}user}, pages pp 273--283, 2000.
\newblock Infinite length modules (Bielefeld, 1998).

\end{thebibliography}
\end{small}

\end{document}